\documentclass[10pt,twoside,reqno]{amsart}
 \usepackage[top=1.4in, bottom=1.4in, left=1in, right=1in]{geometry}
\usepackage{graphicx}
\usepackage{amsmath,latexsym,amssymb}
\usepackage{amsfonts}
\usepackage{mathtools,esint}
\usepackage{amsmath}
\usepackage{amscd,amsthm}
\usepackage{stmaryrd}
\usepackage{fancyhdr}
\usepackage{commath,enumerate}
\usepackage{indentfirst, hyperref}
\usepackage{caption, subcaption}
\usepackage{bbold}

\usepackage{pgf,tikz}
\usepackage{mathrsfs}
\usetikzlibrary{arrows,calc,positioning}
\usepackage{float}
\numberwithin{figure}{section}
\usepackage[bottom]{footmisc}

\numberwithin{equation}{section}

\theoremstyle{definition}

\newcommand{\Rm}[1]{
  \textup{\uppercase\expandafter{\romannumeral#1}}
}

\newcommand{\1}{\mathbb{1}}

\newcommand{\bmo}{\mathrm{BMO}}

\newcommand{\n}{\mathbf{n}}

\newcommand{\R}{\mathbb{R}}

\renewcommand{\u}{\mathbf{u}}
\newcommand{\U}{\mathbf{U}}
\newcommand{\x}{\mathbf{x}}

\renewcommand{\t}{\mathbf{t}}
\newcommand{\y}{\mathbf{x}'}

\def\vp{\varphi}

\def\px{\partial_x}

\newcommand{\diff}{\,\mathrm{d}}

\DeclareMathOperator{\riesz}{\mathbf{\mathcal{R}}}
\DeclareMathOperator{\hilbert}{\mathcal{H}}

\DeclareMathOperator{\sgn}{\mathrm{sgn}}

\def\bel{\begin{equation}\label}
\def\beq{\begin{equation}}
\def\eeq{\end{equation}}
\def\bega{\begin{array}}
\def\enda{\end{array}}

\newcommand{\PV}{\mathrm{PV}}

\newcommand{\Rl}{\mathbb{R}}

\renewcommand{\vec}[1]{\mathbf{#1}}

\author{John K. Hunter}
\address{Department of Mathematics, University of California at Davis}
\email{jkhunter@ucdavis.edu}
\thanks{JKH was supported by the NSF under grant numbers DMS-1616988 and DMS-1908947}
\author{Jingyang Shu}
\address{Department of Mathematics, University of California at Davis}
\email{jyshu@ucdavis.edu}
\author{Qingtian Zhang}
\address{Department of Mathematics, University of California at Davis}
\email{qzhang@math.ucdavis.edu}
\title[Contour Dynamics for SQG Fronts]{Contour Dynamics for Surface Quasi-Geostrophic Fronts}
\date{July 14, 2019}

\begin{document}

\begin{abstract}
We use contour dynamics to derive equations of motion for infinite planar surface quasi-geostrophic (SQG) fronts, and show that it leads to the same result as a regularization procedure introduced previously by Hunter and Shu (2018).
\end{abstract}

\maketitle

\section{Introduction}

In this paper, we use contour dynamics to derive an equation for the
motion of infinite fronts in piecewise constant solutions of the surface quasi-geostrophic (SQG) equation.
The same equation was derived in \cite{HS18} by a regularization procedure
that uses a Galilean transformation to remove a divergence in long-distance  cut-offs of the formal contour dynamics equation.
Thus, the present paper justifies the regularization procedure proposed in \cite{HS18}.
Equations for spatially periodic SQG fronts were also derived by Fefferman and Rodrigo \cite{FR11, Ro05}, and related problems for almost sharp SQG fronts are studied in \cite{CFR04, FLR12, FR12, FR15}.

The SQG equation is a transport equation in two space dimensions $\x = (x, y)$ for an active scalar $\theta(\x, t)$, with
the physical interpretation of a surface buoyancy,
\begin{align}
\label{sqg}
\begin{split}
& \theta_t + \u \cdot \nabla \theta = 0,\qquad
\u = -\riesz^\perp \theta.
\end{split}
\end{align}
The incompressible velocity field $\u(\x, t)$ is determined nonlocally from $\theta(\x,t)$
by a perpendicular Riesz transform $-\riesz^\perp = (\riesz_y,-\riesz_x)$, where $\riesz_x$, $\riesz_y$
are scalar Riesz transforms with respect to $x$, $y$ (see Section \ref{sec:prelim}).
The Riesz transform can also be defined in terms of a Neumann-Dirichlet map for the 3D Laplacian in the derivation of
the 2D SQG equation from the 3D quasi-geostrophic (QG) equation (see Section~\ref{sec:QG}).

The transport equation in \eqref{sqg} preserves piecewise constant solutions in which $\theta(\cdot,t) = \1_{\Omega(t)}(\cdot)$ is the
characteristic function of a domain $\Omega(t) \subset \Rl^2$ with smooth boundary; the boundary moves with normal velocity $\u\cdot\n$ where $\n$ is the normal to the boundary. Contour dynamics, introduced by Zabusky et.~al.~\cite{zabusky}  for the incompressible Euler equations \cite{majda}, allows us to determine the normal velocity from the location of the boundary and derive closed equations for the motion of the boundary.

For the front solutions we consider here, the domain
\begin{equation}
\label{defomega}
\Omega(t) = \{(x,y)\in \Rl^2 : y > \varphi(x,t)\}
\end{equation}
is an upper half-space whose boundary is a graph $y=\varphi(x,t)$, and
\begin{equation}
\theta(\x,t) = \begin{cases} 2\pi & \text{if $y >\varphi(x,t)$},\\0 &\text{if $y<\varphi(x,t)$}.\end{cases}
\label{sqg_front}
\end{equation}
We normalize the jump in $\theta$ across the front to $2\pi$ without loss of generality.
The addition of a constant to $\theta$ does not change the velocity field, so we would get the same result if, for example,
$\theta = \pi \sgn\left[y-\varphi(x,t)\right]$.

As we discuss further in Section~\ref{sec:prelim}, the condition $\u = -\riesz^\perp \theta$ only determines $\u$ up to a spatially uniform flow, and to specify $\u$ uniquely, we require that
\begin{equation}
\u(\x,t) = (2\log |y|,0) + o(1) \qquad \text{as $|y|\to\infty$}.
\label{asyu}
\end{equation}
Our front solutions are then perturbations of the steady SQG shear flow
\begin{equation}
\theta = \begin{cases} 2\pi & \text{if $y >0$},\\0 &\text{if $y<0$},\end{cases}\qquad \u = (2\log |y|,0),
\label{sqg_shear}
\end{equation}
in which disturbances to the flow are caused by the motion of the front and decay away from the front into the interior of the flow.
The expression for $\u$ in  \eqref{sqg_shear} follows from the Hilbert-transform pair \eqref{hilbertsgn}.
The solution \eqref{sqg_shear} is the SQG analog of the linear shear flow $\u = (|y|,0)$ for the 2D incompressible Euler equation with piecewise constant vorticity \cite{BiHu,HS18}.

One can also consider the motion of SQG patches in which $\theta(\cdot,t) = \1_{\Omega(t)}(\cdot)$ is the
characteristic function of a bounded domain $\Omega(t) \subset \Rl^2$, in which case $\theta$ has compact support \cite{CCGS16b, CCG18, Gan08,GP18p,GS14,GS18, KRYZ16, KYZ17}.
An advantage of studying front solutions instead of patches is that they do not introduce extraneous length scales, so they respect the basic scaling properties of the SQG equation and permit an analysis of SQG contour dynamics in a simple geometry. Furthermore, the scalar equation for fronts that can be represented as a graph is simpler than the system of equations for patch-boundaries or fronts that are represented parametrically, although it cannot be used to study front-breaking.  It is reasonable to expect that these front solutions provide an approximation to the motion of sufficiently short wavelength perturbations in patch-boundaries, as well as the local behavior of front-type  solutions in bounded domains which are sufficiently large that the effect of the boundaries on the motion of the front can be neglected.

Unlike the case
of compactly supported patch-solutions for $\theta$, where the far-field velocity can be assumed to approach zero, the far-field velocity of the front solutions is the unbounded flow \eqref{asyu}. The lack of decay in the far-field velocity introduces complications in the reconstruction of the velocity field $\u$ from $\theta$ and the derivation of contour dynamics equations for the front.  The purpose of this paper is to provide
a careful resolution of these complications.

Under suitable assumptions on $\vp$, stated in \eqref{asphi} below, we show that  the location of an SQG front in a solution
\eqref{sqg_front} of  \eqref{sqg} and \eqref{asyu} satisfies
\begin{align}
\label{sqgfront1}
\begin{split}
&\varphi_t(x, t) + 2\left(\log 2 - \gamma\right) \varphi_x(x,t)
\\
&+ \int_\R [\varphi_x(x, t) - \varphi_{x'}(x', t)] \bigg\{\frac{1}{\abs{x-x'}}
- \frac{1}{\sqrt{(x-x')^2 + [\varphi(x, t) - \varphi(x', t)]^2}}\bigg\} \diff{x'} = 2 \log\abs{\partial_x} \varphi_x(x, t),
\end{split}
\end{align}
where $\log|\px|$ is the Fourier multiplier operator with symbol $\log|\xi|$, and $\gamma$ is the Euler-Mascheroni constant. This equation agrees with the one previously derived by a regularization procedure in
 \cite{HS18}, up to a constant $x$-velocity $2\left(\log 2 - \gamma\right) $ which is removed by a Galilean transformation in  \cite{HS18}.
Equation \eqref{sqgfront1} and its generalizations are studied in \cite{HSZ18, HSZ18p, HSZ19pa}.

In this paper, we give two different, but essentially equivalent, derivations of \eqref{sqgfront1}.
The first is based on a decomposition of the solution into a background planar front and a perturbation whose velocity field
approaches zero as $|y|\to \infty$ (see Section \ref{sec:contourI}). The second uses the definition of the Riesz transform
on $L^\infty$-functions to determine the velocity field as the representative in an equivalence class
of $\bmo$-functions with the far-field behavior specified in \eqref{asyu} (see
Section \ref{sec:contourII}).

Before deriving the front equation, we recall some definitions and properties of the Riesz transform in Section~\ref{sec:prelim}, and in Section~\ref{sec:QG}, we explain how it is related to a Neumann-Dirichlet map for the 3D quasi-geostrophic (QG) equation, which is used to derive the 2D SQG equation.

\section{Riesz transforms}
\label{sec:prelim}

In this section we recall some
definitions and properties of the Riesz transform and the space $\bmo$. For more details, see \cite{Duo01, FS72, Ste70, Ste93}.

When $1<p<\infty$, the Riesz transform $\riesz : L^p(\R^n) \to L^p(\R^n; \Rl^n)$ is the bounded singular integral operator defined
pointwise a.e.\ for $f \in L^p(\R^n)$ by (\cite{Duo01}, p.~76)
\begin{align}
\label{rieszdef0}
\begin{split}
\riesz f(\x) &= C_n \mathrm{p.v.}\int_{\R^n} \frac{\x - \y}{|\x - \y|^{n + 1}} f(\y)  \diff{\y}
\\
&= C_n \lim_{\epsilon \to 0+} \int_{\R^n \setminus B_\epsilon(\x)} \frac{\x - \y}{|\x - \y|^{n + 1}} f(\y)  \diff{\y},
\\
C_n &= \frac{1}{\pi^{(n+1)/2}} \Gamma \left(\frac{n+1}{2}\right),
\end{split}
\end{align}
where $B_\epsilon(\x)$ is the ball of radius $\epsilon$ centered at $\x$.   One can also write $\riesz = -\nabla (-\Delta)^{-1/2}$.

For $f \in L^\infty(\R^n)$, the principal value integral on the right-hand side of \eqref{rieszdef0} does not define $\riesz f$, unless it happens to converges absolutely at infinity. However, the Riesz transform can be extended to a bounded linear map $\riesz : L^\infty(\Rl^n) \to \bmo(\Rl^n;\Rl^n)$, where  $\bmo$ denotes the Banach
space of of functions of bounded mean oscillation.

The $\bmo$-norm of $f: \Rl^n\to \Rl$ is defined by
\[
\|f\|_{\bmo} = \sup_{B \subset \R^n} \fint_B \left|f - \fint_B f\right|,\qquad \fint_B f = \frac{1}{|B|} \int_B f,
\]
where $B$ ranges over all balls and $\fint_B f$ denotes the average of $f$ over $B$.
The $\bmo$-norm of a constant is equal to zero, and functions that differ by a constant are regarded as equivalent in $\bmo$.
The space $\bmo$ consists of equivalence classes of locally integrable functions with finite $\bmo$-norms.

The Riesz transform of $f \in L^\infty(\Rl^n)$ can be defined by (\cite{Duo01}, p.~119)
\begin{align}
\label{rieszdef}
\riesz f(\x) = \riesz [f \1_B](\x) + C_n \int_{\R^n \setminus B} \bigg[\frac{\x - \y}{|\x - \y|^{n + 1}} - \frac{\x_0-\y}{|\x_0-\y|^{n + 1}}\bigg]  f(\y) \diff{\y},
\end{align}
where $\x_0\in \Rl^n$ is a fixed point, $B$ is a ball that contains $\x$ and $\x_0$, $\1_B$ is the characteristic function of $B$,
and  $\riesz[f \1_B]$ is defined as in \eqref{rieszdef0}.
The integral on the right-hand side of \eqref{rieszdef} converges absolutely since the integrand is $O(|\y|^{-(n+1)})$ as
$|\y|\to \infty$.
Different choices of $\x_0$ and  $B$ lead to functions that differ by a constant, so they are
equivalent in $\bmo$, and it can be shown that $\riesz f \in \bmo$ for $f \in L^\infty$. In particular, $\riesz 1 = 0$ in $\bmo$. If the support
$E$ of $f \in L^\infty(\Rl^n)$ is a proper subset of $\Rl^n$ and $\x_0 \notin  E$, then we can also define
\begin{align}
\label{rieszdef1}
\riesz f(\x) = C_n \mathrm{p.v.} \int_{E} \bigg[\frac{\x - \y}{|\x - \y|^{n + 1}} - \frac{\x_0-\y}{|\x_0-\y|^{n + 1}}\bigg]  f(\y) \diff{\y},
\end{align}
since this expression agrees with $\eqref{rieszdef}$ up to a constant.

For an SQG shear flow with $\theta = \theta(y)$,
the Riesz transform with respect to $y$ reduces to the Hilbert transform $\hilbert$, and the
corresponding velocity field is $\u = ({u}(y),0)$ where ${u}= \hilbert[{\theta}]$. In particular, if
${\theta}(y)$ is a step function with a jump of $2\pi$, then we have the Hilbert-transform pair (\cite{Duo01}, p.~120)
\begin{align}
\label{hilbertsgn}
{\theta}(y) = \begin{cases} 2\pi &\text{if $y>0$}, \\ 0 &\text{if $y<0$},\end{cases}\qquad
\hilbert[\theta](y) = 2\log|y|,
\end{align}
which gives the planar front solution in \eqref{sqg_shear}.

\section{The quasi-geostrophic equation and Dirichlet-Neumann maps}
\label{sec:QG}

An equivalent way to describe the reconstruction of the SQG velocity field from the buoyancy
is to return to the original derivation of the 2D SQG equation from the 3D QG equation.

In this section, to distinguish between the 2D and 3D variables, we use the notation
\[
\x = (x,y,z),\quad \x_H = (x,y),\quad  \Delta = \partial_x^2 + \partial_y^2 + \partial_z^2, \quad \Delta_H = \partial_x^2 + \partial_y^2,\quad \nabla_H^\perp = (-\partial_y,\partial_x), \quad\riesz_H^\perp = (-\riesz_y,\riesz_x).
\]
The horizontal Riesz transform $-\riesz_H^\perp = \nabla_H^\perp (-\Delta_H)^{-1/2}$
in the SQG equation then arises as the orthogonal gradient of a Neumann-Dirichet map $(-\Delta_H)^{-1/2}$
for the 3D Laplacian in the QG equation.
The particular choice for  the Riesz transform is determined by the far-field boundary conditions for the QG equation.

The QG equation provides an approximate description of nearly horizontal geostrophic flows in a vertically stratified fluid \cite{majda_qg, Ped87}.
In suitably non-dimensionalized variables, the streamfunction $\Psi(\x,t)$ of the flow satisfies $\Delta \Psi = \PV$, where $\PV$ is the potential
vorticity in the fluid. The horizontal velocity of the fluid is $\U_H = \nabla_H^\perp \Psi$. The streamfunction is proportional to the fluid pressure, and $\Psi_z$ has the interpretation of a temperature perturbation or buoyancy, rather than a vertical velocity component.

The SQG equation describes quasi-geostrophic flows in a half-space $\Rl^2\times \Rl^+$ with zero potential vorticity
in $z>0$ and a temperature jump, or surface buoyancy, $\theta(\x_H,t)$ at $z=0$, which is transported by the velocity field
$\u_H = \left.\U_H\right|_{z=0}$ on the boundary \cite{sqg_lap,Ped87}:
\begin{align*}
\begin{split}
&\Delta \Psi = 0\quad \text{in $z>0$}, \qquad -\left.\partial_z \Psi\right|_{z=0} = \theta,
\\
&\theta_t + \u_H\cdot \nabla_H \theta = 0,\qquad \u_H = \nabla_H^\perp \left.\Psi\right|_{z=0}.
\end{split}
\end{align*}

We omit an explicit indication of the time-variable, and write $\psi = \left.\Psi\right|_{z=0}$ for the value of $\Psi$ on the boundary.
Then $\u_H = \nabla_H \psi$ and $\psi$ is related to $\theta$ by
a solution of the Neumann problem
\begin{align}
\begin{split}
&\Delta \Psi = 0\quad \text{in $z>0$},
\qquad
-\left.\partial_z \Psi\right|_{z=0} = \theta,
\end{split}
\label{QG_DN}
\end{align}
meaning that $\theta\mapsto \psi$ is a Neumann-Dirichlet map for the 3D Laplacian in the upper half space.
From the point of view of potential theory, this problem is the same as finding the electrostatic potential $\Psi$
of a semi-infinite charged plate located at $\{(x,y,0) \in \Rl^3 : y > \vp(x)\}$ with a constant surface charge density of $4\pi$.

The solution of \eqref{QG_DN} is unique up to a harmonic function $\Psi'(\x)$ in $z>0$ with zero normal derivative on $z=0$, which can be fixed by imposing suitable boundary conditions at infinity.
For example, the addition of a linear harmonic function $\Psi' = Ax + By$ to $\Psi$ does not change $\theta$ and adds a uniform velocity field $\u_H = (-B,A)$ to $\u_H$. On the other hand, if $\theta = C$ is constant, then the solution $\Psi' = Cz + D$ (corresponding to a uniform temperature gradient in the QG equation) gives $\psi = 0$, so the addition of a constant
to $\theta$ has no effect on the corresponding velocity field $\u_H$.

In particular, let us consider the QG solution that corresponds to the planar front solution in \eqref{sqg_shear}, with
\begin{equation}
\theta(\x) = \begin{cases} 2\pi & \text{if $y >0$},\\0 &\text{if $y<0$}.\end{cases}
\label{thetabak}
\end{equation}
Differentiating \eqref{QG_DN} with respect to $z$, we see that ${\Phi} = {\Psi_z}$ satisfies the Dirichlet problem
\begin{align}
\begin{split}
&\Delta {\Phi} = 0 \quad \text{in $z > 0$},\qquad
- {\Phi} \Big|_{z = 0} = {\theta}.
\end{split}
\label{PhiBVP}
\end{align}
We look for solutions of \eqref{thetabak}--\eqref{PhiBVP} that are independent of $x$. Then $\Phi_{yy} + \Phi_{zz} = 0$, whose
general solution for the Fourier transform of $\Phi(y,z)$ with respect to $y$,
\[
\hat{\Phi}(\xi,z) = \frac{1}{2 \pi} \int_{\R} \Phi(y,z) e^{-i \xi y}\diff{y},
\]
is given by $\hat{\Phi}(\xi, z) = A(\xi) e^{-|\xi| z} + B(\xi) e^{|\xi| z}$.

We further require that $\hat{\Phi}(\xi, z)\to 0$ as $z \to \infty$ for $\xi\ne 0$, in which case $B = 0$ and $\hat{\Phi}(\xi, z) = \hat{\theta}(\xi) e^{-|\xi| z}$. Inverting this Fourier transform, we get that
\begin{align*}
{\Phi}(y, z) = \pi + 2 \arctan\bigg(\frac{y}{z}\bigg),
\end{align*}
and taking an antiderivative of ${\Phi}$ with respect to $z$, we get the streamfunction
\[
{\Psi}(y, z) = -2y + y \log\left(y^2 + z^2\right) + 2 z \arctan\left(\frac{y}{z}\right) + \pi z.
\]
This function provides the appropriate far-field behavior as $y^2+z^2 \to \infty$ of QG-front solutions
in defining the Neumann-Dirichlet map from \eqref{QG_DN}.

The boundary value of ${\Psi}$ on $z = 0$ is
\[
{\psi}(y) = \lim_{z \to 0^+} {\Psi}(y, z) = -2 y + 2 y \log|y|,
\]
with the velocity field ${\u}_H  = \left(2 \log|y|, 0\right)$,  as in the planar front solution \eqref{sqg_shear}.

\section{Contour dynamics equation I}
\label{sec:contourI}

We now derive contour dynamics equations for the front solutions \eqref{defomega}--\eqref{sqg_front} by decomposing the solution
into a planar shear flow and a perturbation whose velocity field approaches zero as $|y|\to\infty$.

We denote the front by $\Gamma(t) = \partial\Omega(t)$, and consider its motion on a
time interval $0\le t \le T$ for some $T>0$. We assume that:
\begin{align}
\begin{split}
&\text{(i) $\varphi(\cdot,t) \in C^{1,\alpha}(\Rl)$ for some $\alpha > 0$ and $\varphi(x,t)$ is bounded
on $\Rl\times [0,T]$};
\\
&\text{(ii) $\varphi_x(x,t) = O(|x|^{-\beta})$ as $|x|\to \infty$ for some $\beta > 0$}.
\end{split}
\label{asphi}
\end{align}
In that case, all of the integrals in the following converge.

We choose $h > 0$ such that $- h < \inf \{\vp(x,t) : (x,t) \in \R\times[0,T]\}$, and let
\begin{align}
\label{thetatilde}
\tilde{\theta}(\x) = \begin{cases}2 \pi & \text{if $\ y > -h$},\\ 0 & \text{if $y < -h$},\end{cases},
\qquad
\tilde{\u}(\x) = \left(2 \log|y + h|, 0\right),
\end{align}
be the planar front solution \eqref{sqg_shear} translated to $y=-h$.

\begin{figure}[ht]
\centering
\includegraphics[width=0.7\textwidth]{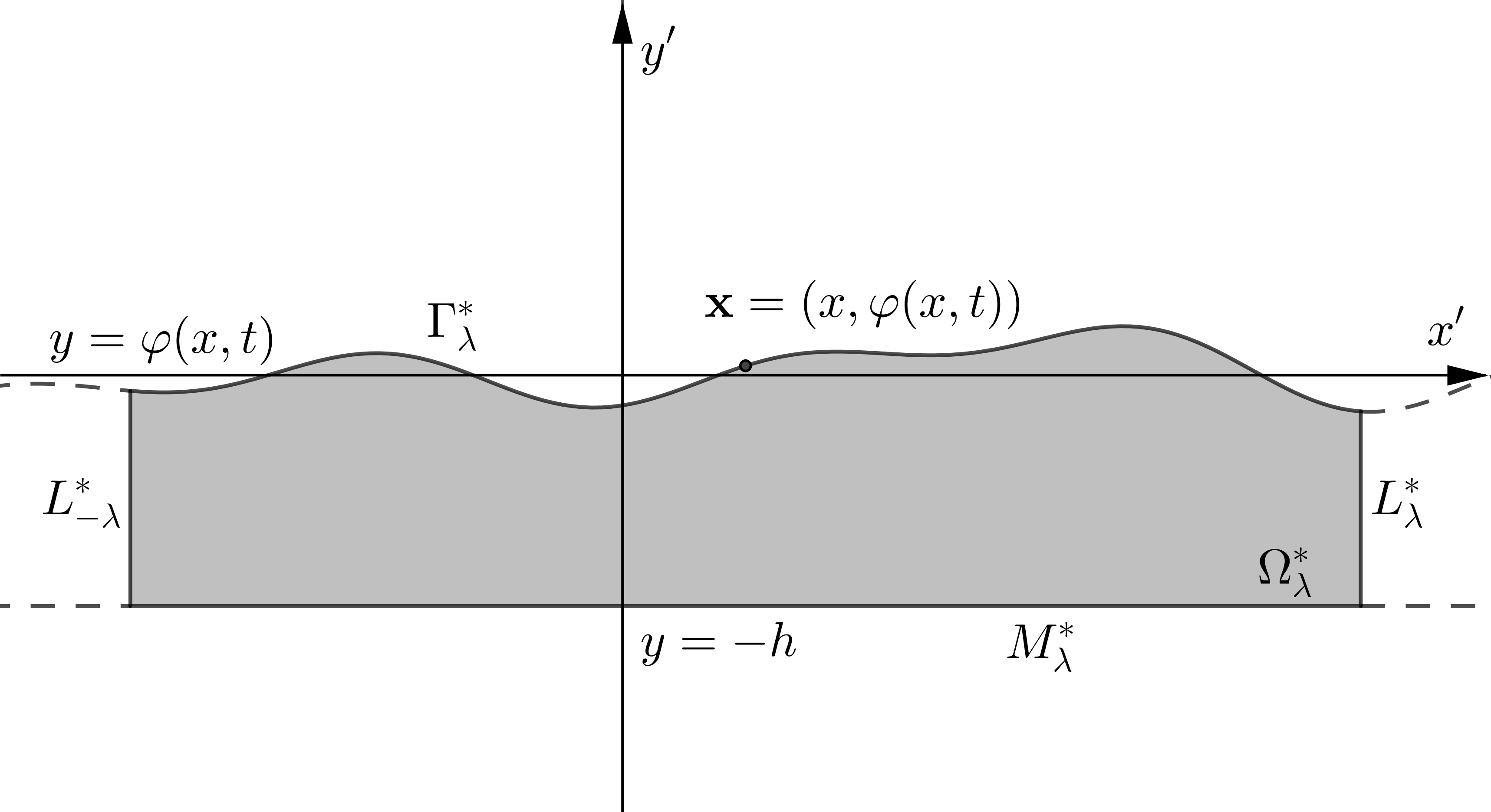}
\caption{An illustration of the cut-off region $\Omega_\lambda^*$ in \eqref{defomega*} with a point $\x$ on the front. The boundary $\partial\Omega^*_\lambda$ consists of the lines $L^*_{\pm\lambda}: x' = x\pm\lambda$ with $-h \leq y' \leq \varphi(x\pm\lambda)$,
$M^*_\lambda: y' = -h$ with $|x-x'| \le \lambda$, and the cut-off front $\Gamma^*_\lambda : y'=\vp(x')$ with $|x-x'|\le \lambda$.
The function $\theta^*$ in \eqref{deftheta*}  is equal to $-2 \pi$ in the strip $-h < y < \vp(x, t)$) and equal to $0$ in $y > \vp(x, t)$ or $y < -h$.}
\label{fig:fronts_diff2}
\end{figure}

We decompose the front solution \eqref{sqg_front} as
\begin{align*}
\theta(\x,t) = \tilde{\theta}(\x) + \theta^*(\x,t),
\end{align*}
where $\tilde{\theta}$ is defined in \eqref{thetatilde}, and
\begin{equation}
 \theta^*(\x,t) = \begin{cases} -2\pi &\text{if $-h < y < \vp(x,t)$},\\
 0 &\text{otherwise}.\end{cases}
 \label{deftheta*}
 \end{equation}
We denote the support of $\theta^*(\cdot,t)$ by $\Omega^*(t)$.
The corresponding decomposition of the velocity field is
\begin{equation}
\u(\x,t) = \tilde{\u}(\x) + \u^*(\x,t),
\label{decu}
\end{equation}
where $\tilde{\u}$ is defined in \eqref{thetatilde} and $\u^* = -\riesz^\perp \theta^*$ is given by
\[
\u^*(\x,t) = \mathrm{p.v.} \int_{\Omega^*(t)} \frac{(\x - \x')^\perp}{|\x - \x'|^3} \diff{\x'},\qquad (x,y)^\perp = (-y,x).
\]
This principal value integral converges absolutely at infinity, since, writing $\x' = (x', y')$, the integrand is $O\left(|x'|^{-2}\right)$ as
$|x'|\to \infty$ and compactly supported in $y'$. It follows that
\begin{equation}
\u^*(\x) = \lim_{\lambda \to \infty} \u^*_\lambda(\x), \qquad \u^*_\lambda(\x) = \mathrm{p.v.} \int_{\Omega^*_\lambda(x,t)} \frac{(\x - \x')^\perp}{|\x - \x'|^3} \diff{\x'},
\label{defulambda}
\end{equation}
where (see Figure \ref{fig:fronts_diff2})
\begin{equation}
\Omega^*_\lambda(x,t) = \left\{\x' \in \R^2 : \text{$|x - x'| < \lambda$, $-h < y' < \vp(x',t)$}\right\}.
\label{defomega*}
\end{equation}

Let $\x = \left(x,\varphi(x)\right)$ be a point on the front and denote by
\begin{equation}
\vec{n}(\x,t) = \frac{1}{\sqrt{1+\varphi_x^2(x,t)}}(-\varphi_x(x,t),1)
\label{defn}
\end{equation}
the unit upward normal to $\Gamma(t)$ at $\x$. The motion of the front is determined by the normal
velocity $\u\cdot \n$, which is continuous and well-defined on the front.
The tangential component of $\u$ diverges to infinity,
but this does not affect the motion of the front.

We take the inner product of $\u^*_\lambda$ in \eqref{defulambda} with $\n$, write
\[
\frac{(\x - \y)^\perp}{|\x - \y|^{3}}
= \nabla_{\y}^\perp  \frac{1}{|\x-\y|},
\]
and apply Green's theorem, to get that
\begin{align}
\label{ulambda*}
\u^*_\lambda(\x,t)\cdot\n(\x,t) = - \int_{\partial \Omega^*_\lambda(x,t)} \frac{\t(\x',t)\cdot \n(\x,t)}{|\x - \x'|} \diff{s(\x')},
\end{align}
where $\t$ is the negatively oriented unit tangent vector on $\partial \Omega^*_\lambda$ and $\diff{s(\x')}$ is an element of arclength.
Since $\t(\x,t)\cdot \n(\x,t) = 0$, the assumed H\"older continuity of $\varphi_x$ ensures that this integral converges
at $\y = \x$, so there is no contribution from the principal value at $\x' = (x,\vp(x,t))$.

As illustrated in Figure \ref{fig:fronts_diff2}, we decompose the boundary
as $\partial \Omega^*_\lambda = \Gamma^*_\lambda \cup M_\lambda^* \cup L_{\lambda}^*\cup L^*_{-\lambda}$.
On $L_{-\lambda}^*$, we have $\t(\x',t)=(0,1)$, $x'=-\lambda$, and $\diff{s(\x')}=dy'$, so
\[
\int_{L_{-\lambda}^*} \frac{\t(\x',t)}{|\x - \x'|} \diff{s(\x')}  = (0,1) I^*_\lambda(x,t),
\]
where
\begin{align*}
I^*_\lambda(x,t) &=\int_{-h}^{\vp(-\lambda)} \frac{1}{\sqrt{(x + \lambda)^2 + (\vp(x,t) - y')^2}} \diff{y'}
\\
&= \log\left(\vp(x,t) - \vp(-\lambda,t) + \sqrt{(x + \lambda)^2 + (\vp(x,t) - \vp(-\lambda,t))^2}\right)\\
& \qquad - \log\left(\vp(x,t) + h + \sqrt{(x + \lambda)^2 + (\vp(x,t) + h)^2}\right)\\
& \to 0 \qquad \text{as $\lambda \to \infty$},
\end{align*}
since $\vp$ is bounded. Similarly, the limit of the integral over $L^*_\lambda$ as $\lambda\to \infty$ also vanishes, so the only contribution
to $\u^*$ comes from $\Gamma^*_\lambda$ and $M^*_\lambda$.

The tangent vector on $\Gamma^*_\lambda$ is
\begin{equation}
\vec{t}(\y,t) = \frac{1}{\sqrt{1+\varphi_{x'}^2(x',t)}}\left(1,\varphi_{x'}(x',t)\right),
\label{deft}
\end{equation}
and the tangent vector on $M^*_\lambda$ is $(-1,0)$. Using \eqref{defn} and \eqref{deft} in \eqref{ulambda*}, and taking the limit $\lambda\to\infty$, we get that
\begin{align*}
\u^*(\x, t) \cdot \n(\x, t) &= - \lim_{\lambda\to \infty} \int_{\Gamma^*_\lambda \cup M^*_\lambda}
\frac{\t(\x',t)\cdot \n(\x,t)}{|\x - \x'|} \diff{s(\x')}
=\frac{1}{\sqrt{1 + \vp_x^2(x, t)}} I^*(x, t),\\
I^*(x, t) &= \int_\R \left\{\frac{\vp_x(x, t) - \vp_{x'}(x', t)}{\sqrt{(x - x')^2 + (\vp(x, t) - \vp(x', t))^2}}
- \frac{\vp_x(x, t)}{\sqrt{(x - x')^2 + (\vp(x, t) + h)^2}}\right\} \diff{x'}.
\end{align*}

Including the contribution from the background flow $\tilde{\u}$,
and using the condition that the front $y = \vp(x, t)$ moves with the upward normal velocity $\u \cdot \n = (\tilde{\u} + \u^*) \cdot \n$,
we obtain that
\[
\vp_t(x, t) = I(x, t),\qquad I(x, t) = \left(\sqrt{1+ \vp_x^2}\right)\tilde{\u}(\x, t) \cdot \n(\x, t) + I^*(x, t).
\]
From \eqref{thetatilde} and \eqref{defn}, we have
\[
 \left(\sqrt{1+ \vp_x^2}\right) \tilde{\u}(\x, t) \cdot \n(\x, t)  = - 2 \left(\log|\vp(x, t) + h|\right) \vp_x(x, t).
\]
We then decompose $I$ as
\begin{align}
\begin{split}
I(x,t) &= I_1(x, t) + I_2(x, t) + I_3(x, t),
\\
I_1(x, t) & = \int_\R \left\{\frac{\vp_x(x, t) - \vp_x(x', t)}{\sqrt{(x - x')^2 + (\vp(x, t) - \vp(x', t))^2}}
- \frac{\vp_x(x, t) - \vp_x(x', t)}{|x - x'|}\right\} \diff{x'},\\
I_2(x, t) & = \int_\R \left\{\frac{\vp_x(x, t) - \vp_x(x', t)}{|x - x'|} - \frac{\vp_x(x, t)}{\sqrt{(x')^2 + 1}}\right\} \diff{x'},\\
I_3(x, t) & = \vp_x(x, t) \left\{\int_\R \frac{1}{\sqrt{(x')^2 + 1}}
- \frac{1}{\sqrt{(x - x')^2 + (\vp(x, t) + h)^2}}\diff{x'} - 2 \log|\vp(x, t) + h|\right\}.
\end{split}
\label{defI}
\end{align}
All of these integrals converge in view of the assumed H\"older continuity and decay of $\vp_x$.

Direct evaluation of the integral for $I_3$ yields
\begin{align*}
I_3(x, t) &= \vp_x(x, t) \left\{\Big[\sinh^{-1}{x'} - \log\left(x' - x + \sqrt{(x' - x)^2 + (\vp(x, t) + h)^2}\right)\Big]_{x' = -\infty}^\infty
- 2 \log|\vp(x, t) + h|\right\}
\\
&= 0.
\end{align*}
Thus, the equation for the front is
\begin{equation}
\vp_t(x, t) = I_1(x, t) + I_2(x, t),
\label{fronteq}
\end{equation}
where $I_1$ is the nonlinear term and $I_2$ is the linear term in the equation.

To express $I_2$ as a Fourier multiplier operator, we note that $I_2=0$ if $\vp =1$, and
if $\vp(x) = e^{i\xi x}$ with $\xi\ne 0$, then
\begin{align*}
I_2(x) &= i\xi e^{i\xi x} \int_\R \left[ \frac{ 1- e^{i\xi(x'-x)}}{|x-x'|} - \frac{ 1}{\sqrt{(x')^2 + 1}} \right] \diff{x'}
\\
&= i\xi e^{i\xi x}\int_\Rl\left[ \frac{ 1- e^{i\xi(x'-x)}}{|x-x'|} - \frac{ 1}{\sqrt{(x-x')^2 + 1}} \right] \diff{x'}
\\
& = 2i\xi e^{i\xi x}\int_0^\infty\left[ \frac{ 1- \cos \xi s}{s} - \frac{ 1}{\sqrt{s^2 + 1}} \right] \diff{s}.
\end{align*}
Using the identity
\begin{equation}
\int_0^\infty\left[\frac{ 1}{\sqrt{s^2 + 1}} - \frac{1}{\sqrt{s^2 + c^2}} \right] \diff{s} = \log|c|,
\label{scaleid}
\end{equation}
with $c= 1/|\xi|$, the change of variable $s' = |\xi|s$, and a cosine integral, we get that
\begin{align*}
I_2(x)
& = 2i\xi e^{i\xi x}\left(\log|\xi| + \int_0^\infty\left[ \frac{ 1- \cos s'}{s'} - \frac{ 1}{\sqrt{(s')^2 + 1}} \right] \diff{s'}\right)
\\
&=2i\xi e^{i\xi x}\left(\log|\xi| + \gamma -\log 2\right),
\end{align*}
where $\gamma$ is the Euler-Mascheroni constant. It follows that
\[
I_2 = 2\log|\partial_x| \vp_x + 2(\gamma -\log 2)\vp_x,
\]
where $ \log|\partial_x|$ is the Fourier multiplier operator with symbol
$\log|\xi|$.
Thus, using the expressions for $I_1$, $I_2$ in \eqref{fronteq}, we get the front equation \eqref{sqgfront1}.

As in Section~\ref{sec:contourII}, one can verify that $\u^*(\x,t) = o(1)$ as $|y|\to \infty$. Since
\[
2\log |y+h| = 2\log|y| + o(1)\qquad \text{as $|y|\to \infty$},
\]
the velocity field $\u(\x,t)$ in \eqref{decu} corresponding to this solution satisfies \eqref{asyu}.

The appearance of a logarithm in the far-field boundary condition \eqref{asyu} breaks the scale-invariance of the SQG equation under $\x \mapsto k \x$, $t \mapsto k t$ for $k > 0$. Instead, the invariance that preserves the scale-invariance of the SQG equation \eqref{sqg} and the asymptotic behavior
\eqref{asyu} of the velocity, is a combined scaling-Galilean transformation \cite{HS18}
\[
x\mapsto k [x + (2\log|k|)t],\qquad y \mapsto ky,\qquad t \mapsto kt.
\]
This behavior is  consistent with the scale-invariance of the Riesz transform on
$L^\infty$, where the rescaling of an $L^\infty$-function (such as $\pi \sgn y$) may map a representative of the Riesz transform (such as $2\log|y|$) to an equivalent representative in $\bmo$ that differs from the original one by a constant.

From the point of view of dimensional analysis, both $\theta$ and $\u$ in the SQG equation have the dimension of velocity,
the Riesz transform being a dimensionless operator. The only dimensional parameter in the front data is a velocity,
namely the jump in $\theta$ across the front, which is non-dimensionalized to $2\pi$ in \eqref{sqg_front}.
There are no intrinsic length or time scales, but the spatial variables are implicitly
non-dimensionalized by the condition that the asymptotic far-field velocity in \eqref{asyu} vanishes at $y=1$, a condition that is
scale-Galilean invariant, but not scale invariant.

Similar issues are well-known in potential theory for unbounded charge distributions.
For example, there is no length scale in the problem for the
electrostatic potential of an infinite charged wire, which is given by the logarithmic Newtonian potential.
The potential diverges at infinity, so one cannot normalize a zero point for the potential by requiring
that the potential approaches zero at infinity (as one usually does for compact charge distributions).
Instead, one picks an arbitrary radial distance $r_0 > 0$ from the wire and requires that the potential
vanish at a distance $r=r_0$, or $r=1$ in spatial variables non-dimensionalized by $r_0$ (see e.g.\ Sec.~III.5 in Kellogg \cite{kellogg}).
The problem is then invariant under spatial rescalings and an appropriate shift in the zero-point of the potential.

\section{Contour dynamics equation II}
\label{sec:contourII}

In this section, we give an alternative derivation of \eqref{sqgfront1} based on the definition of the $\bmo$-valued Riesz transform on $L^\infty$
in \eqref{rieszdef1}.
As before, we assume that $\vp$ satisfies \eqref{asphi}.

From \eqref{sqg} and \eqref{rieszdef1}, with $n=2$ and $C_2= 1/2\pi$, a representative velocity field of the front solution \eqref{sqg_front} is given by
\begin{align}
\u(\x,t) =  -\text{p.v.} \int_{\Omega(t)} \left[\frac{(\x - \y)^\perp}{|\x - \y|^{3}} - \frac{(\x_0-\y)^\perp}{|\x_0-\y|^{3}}\right] \diff{\y}
- \bar{\u}(t),
\label{defuu}
\end{align}
where $\Omega(t)$ is given by \eqref{defomega} and $\x_0\notin \bar{\Omega}(t)$.
For definiteness, we choose $\x_0 = (0,-h)$
where $h > 0$ and  $-h < \inf \{\varphi(x,t) : (x,t)\in \Rl\times[0,T]\}$.
The spatially uniform velocity $\bar{\u}(t) = (\bar{u}(t),\bar{v}(t))$ in \eqref{defuu}
will be chosen so that $\u(\x,t)$ satisfies the far-field condition \eqref{asyu}.
However, any such representative leads to equivalent dynamics for the fronts, since any
uniform velocity $\left(\bar{u}(t),\bar{v}(t)\right)$
can be removed by a translation $(x,y) \mapsto (x-a(t), y-b(t))$ where $(a,b)_t = (\bar{u}, \bar{v})$.

\begin{figure}[h]
\centering
\includegraphics[width=0.7\textwidth]{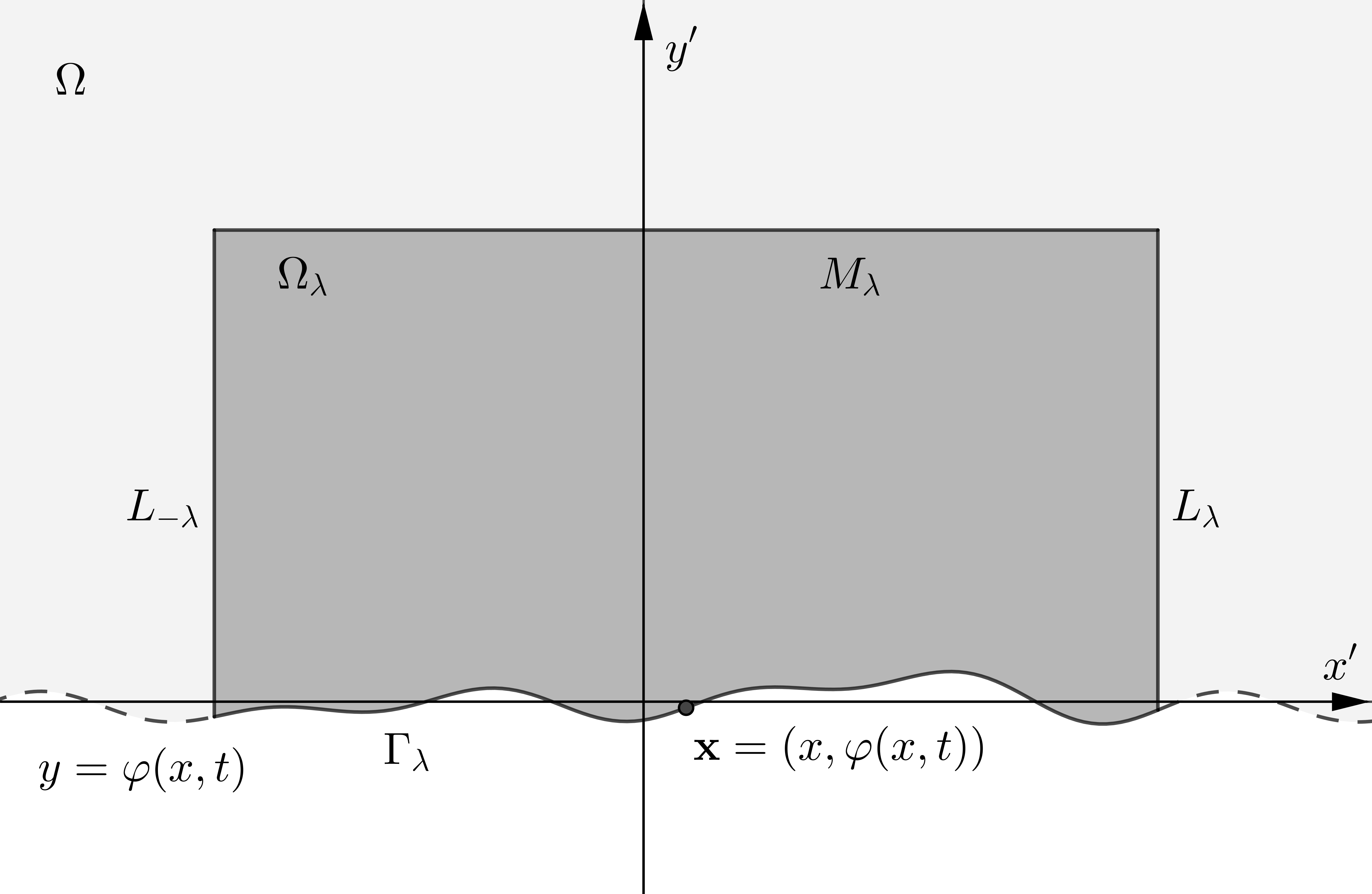}
\caption{An illustration of the cut-off region $\Omega_\lambda$ in \eqref{defomegal} with a point $\x$ on the front. The boundary $\partial\Omega_\lambda$ consists of the lines $L_{\pm\lambda}: x' = x\pm\lambda$ with $\varphi(x\pm\lambda)\le y' \le \lambda$,
$M_\lambda: y' = \lambda$ with $|x-x'| \le \lambda$, and the cut-off front $\Gamma_\lambda : y'=\vp(x')$ with $|x-x'|\le \lambda$.}
\label{fig:lambda_cutoff}
\end{figure}

Since the integral in \eqref{defuu} converges absolutely at infinity, we have
 \begin{equation*}
 \u(\x,t)  = \lim_{\lambda\to \infty}  \u_\lambda(\x,t) - \bar{\u}(t),\qquad
\u_\lambda(\x,t) = -\text{p.v.}\int_{\Omega_\lambda(x,t)} \left[\frac{(\x - \y)^\perp}{|\x - \y|^{3}} - \frac{(\x_0-\y)^\perp}{|\x_0-\y|^{3}}\right] \diff{\y},
\end{equation*}
where (See Figure~\ref{fig:lambda_cutoff})
\begin{equation}
\Omega_\lambda(x,t)  = \left\{\x'\in \Rl^2 : \text{$|x'-x| < \lambda$, $\varphi(x',t) < y' <\lambda$}\right\}.
\label{defomegal}
\end{equation}

First, we consider the case when $\x\notin \Gamma(t)$. We write
\[
\frac{(\x - \y)^\perp}{|\x - \y|^{3}} - \frac{(\x_0-\y)^\perp}{|\x_0-\y|^{3}}
= \nabla_{\y}^\perp \left[ \frac{1}{|\x-\y|} - \frac{1}{|\x_0-\y|}\right],
\]
and apply Green's theorem to get that
\begin{equation}
\u_\lambda(\x,t) = - \int_{\partial \Omega_{\lambda}(x,t)}\left[\frac{1}{|\x-\y|} - \frac{1}{|\x_0-\y|}\right]  \t(\y,t) \diff{s(\y)},
\label{ulambdaeq}
\end{equation}
where $\t(\y,t)$ is the positively
oriented unit tangent vector on $\partial \Omega_{\lambda}(x,t)$ and $\diff{s(\x')}$ is an element of arclength. There is no contribution from the
principal value, since the corresponding integral of $\t(\y)/|\x-\y|$ over $\partial B_\epsilon(\x)$ is zero.

As illustrated in Figure~\ref{fig:lambda_cutoff}, we decompose the boundary as $\partial \Omega_\lambda(x,t) = \Gamma_\lambda(x,t) \cup C_\lambda(x,t)$, where  $C_\lambda(x,t)$ consist of the lines $L_{\pm\lambda}(x,t)$ and $M_\lambda(x,t)$, and $\Gamma_\lambda(x,t)$ is the cut-off front
with tangent vector \eqref{deft}.
The integrand in \eqref{ulambdaeq} is $O(\lambda^{-2})$ on $C_\lambda$, so taking the limit of  \eqref{ulambdaeq} as $\lambda \to \infty$, we get for $\x\notin \Gamma$ that
\begin{equation}
\u(\x,t) = -\int_{\Gamma(t)}\left[\frac{1}{|\x-\y|} - \frac{1}{|\x_0-\y|}\right]  \t(\y,t) \diff{s(\y)} - \bar{\u}(t).
\label{ueq1}
\end{equation}

Writing out \eqref{ueq1} in components, we find that
\begin{align}
u(x,y,t) &= -\int_{\Rl} \left\{\frac{1}{\sqrt{(x-x')^2 + (y-\varphi(x',t))^2}} - \frac{1}{\sqrt{(x-x')^2 + (h+\varphi(x',t))^2}}\right\} \, \diff{x'} - \bar{u}(t),
\label{defu}
\\
v(x,y,t) &= -\int_{\Rl} \left\{\frac{1}{\sqrt{(x-x')^2 + (y-\varphi(x',t))^2}} - \frac{1}{\sqrt{(x-x')^2 + (h+\varphi(x',t))^2}}\right\} \varphi_{x'}(x',t)\, \diff{x'} - \bar{v}(t).
\label{defv}
\end{align}
From \eqref{defu}, we see that
\begin{align*}
u(x,y,t) = - \int_{\Rl} \left\{\frac{1}{\sqrt{(x')^2 + y^2}} - \frac{1}{\sqrt{(x')^2 + (h + \varphi(x',t))^2}} \right\} \, \diff{x'} - \bar{u}(t) + o(1)
\qquad\text{as $|y|\to\infty$}.
\end{align*}
Using the identity \eqref{scaleid}
with $c = |y|$, we get that
\begin{align*}
&-\int_{\Rl} \left\{ \frac{1}{\sqrt{(x')^2 + y^2}} - \frac{1}{\sqrt{(x')^2 + (h + \varphi(x'))^2}} \right\} \, \diff{x'}
 \\&\qquad = 2\log|y| - \int_{\Rl} \left\{\frac{1}{\sqrt{(x')^2 + 1}} - \frac{1}{\sqrt{(x')^2 + (h + \varphi(x'))^2}}\right\} \, \diff{x'}.
\end{align*}
It follows that $u (x,y,t) = 2\log|y| + o(1)$ as $|y|\to \infty$ if
\begin{equation}
\bar{u}(t) = -\int_\R \left[ \frac{1}{\sqrt{(x')^2 + 1}} - \frac{1}{\sqrt{(x')^2 + (h + \vp(x',t))^2}} \right] \diff{x'}.
\label{defubar}
\end{equation}

We also see directly from \eqref{defv} that
\[
v(x,y,t) =  \int_\R  \frac{\vp_{x'}(x',t)}{\sqrt{(x')^2 + (h + \vp(x',t))^2}}\diff{x'}
-\bar{v}(t) + o(1) \qquad\text{as $|y|\to \infty$},
\]
so $v (x,y,t) = o(1)$ as $|y|\to \infty$ if
\begin{equation}
\bar{v}(t) =  \int_\R  \frac{\vp_{x'}(x',t)}{\sqrt{(x')^2 + (h + \vp(x',t))^2}}\diff{x'}.
\label{defvbar}
\end{equation}
The velocity field \eqref{ueq1} therefore satisfies \eqref{asyu} when $\bar{\u} = (\bar{u},\bar{v})$ is given by
\eqref{defubar}--\eqref{defvbar}.

Next, let $\x = \left(x,\varphi(x)\right)$ be a point on the front $\Gamma(t)$, with upward normal \eqref{defn}.
We take the inner product of $\u_\lambda$ in \eqref{ulambdaeq} with $\n$ and take the limit $\lambda\to\infty$ as before,  to get that
\begin{equation*}
\u(\x,t)\cdot \n(\x,t) = - \int_{\Gamma(t)}\left[\frac{1}{|\x-\y|} - \frac{1}{|\x_0-\y|}\right]  \t(\y,t)\cdot\n(\x,t) \diff{s(\y)} - \bar{\u}(t)\cdot \n(\x)
\end{equation*}
Writing out this integral explicitly and using the expressions for $\t$, $\n$, $\x$, and $\x_0$,
we find that the normal velocity on the front is
\begin{align}
\begin{split}
&\u(\x,t)\cdot\n(\x,t) = \frac{1}{\sqrt{1 + \vp_x^2(x,t)}} \left[J(x,t) + \vp_x(x,t) \bar{u}(t) - \bar{v}(t)\right],
\\
&J(x,t) =  \int_\R \left[\frac{ \vp_x(x,t)-\vp_{x'}(x',t)}{\sqrt{(x - x')^2 + (\vp(x,t) - \vp(x',t))^2}} - \frac{\vp_{x}(x,t) - \vp_{x'}(x',t)}{\sqrt{(x')^2 + (h+\vp(x',t))^2}}\right]\diff{x'}.
\end{split}
\label{un1}
\end{align}
The condition that the front $y=\vp(x,t)$ moves with the normal velocity $\u\cdot\n$ implies that
\begin{equation}
\vp_t(x,t) = J(x,t)+ \vp_x(x,t) \bar{u}(t) - \bar{v}(t).
\label{un_front1}
\end{equation}

We decompose the integral for $J$ in \eqref{un1} as
\begin{align*}
J(x,t)  = I_1(x,t) + I_2(x,t) + J_3(x,t) + J_4(t),
\end{align*}
where $I_1$, $I_2$ are given in \eqref{defI} and
\begin{align*}
J_3(x,t) &= \vp_x(x,t)\int_\R \left[ \frac{1}{\sqrt{(x')^2 + 1}} - \frac{1}{\sqrt{(x')^2 + (h + \vp(x',t))^2}} \right] \diff{x'},
\\
J_4(t)  &= \int_\R  \frac{\vp_{x'}(x',t)}{\sqrt{(x')^2 + (h + \vp(x',t))^2}}\diff{x'}.
\end{align*}
From \eqref{defubar}--\eqref{defvbar}, we see that $J_3 = -\vp_x\bar{u}$ and $J_4 = \bar{v}$, so \eqref{un_front1} becomes
\eqref{fronteq}, and we get the same result as before.

\end{document}